\documentclass{amsart}

\usepackage{latexsym,enumerate}
\usepackage{amsmath,amsthm,amsopn,amstext,amscd,amsfonts,amssymb}
\usepackage[ansinew]{inputenc}
\usepackage{verbatim}
\usepackage{epsfig} 
\usepackage{graphicx,psfrag} 
\usepackage{subfigure}

\newcommand{\CC}{\mathbb{C}}

\newcommand{\RR}{\mathbb{R}}
\newcommand{\ZZ}{\mathbb{Z}}

\renewcommand{\SS}{\mathbb{S}}

\newcommand{\Ep}{\mathcal{E}}

\newtheorem{teorema}{Theorem}[section]
\newtheorem{lema}[teorema]{Lemma}
\newtheorem{prop}[teorema]{Proposition}
\newtheorem{corolario}[teorema]{Corollary}
\newtheorem{definicion}[teorema]{Definition}

\newtheorem{obs}[teorema]{Remark}

\DeclareMathOperator{\dist}{dist}

\newcommand{\spb}[1]{\smallskip}
\newcommand{\mpb}[1]{\medskip}
\newcommand{\bpb}[1]{\bigskip}

\newcommand{\p}{\partial}

\renewcommand{\a}{\alpha}
\renewcommand{\b}{\beta}
\newcommand{\e}{\varepsilon}
\renewcommand{\d}{\delta}

\newcommand{\g}{\gamma}
\renewcommand{\th}{\theta}

\newcommand{\s}{\sigma}

\begin{document}

\title{STRUCTURE THEOREM FOR RIEMANNIAN SURFACES
WITH ARBITRARY CURVATURE}

\author{
Ana Portilla$^{(1)(2)(3)}$, Jose M. Rodriguez
$^{(1)(2)(3)}$ \and Eva Touris$^{(1)(2)(3)}$}
\date{February 28, 2007.
$(1)\,\,\,$ Research partially supported by a grant from
DGI (BFM 2006-11976), Spain.\\
$(2)\,\,\,$ Research partially supported by a grant from
DGI (BFM 2006-13000-C03-02), Spain. \\
$(3)\,\,\,$ Research partially supported by a grant from
MEC (MTM 2006-26627-E), Spain.}
\email{apferrei@math.uc3m.es, jomaro@math.uc3m.es, etouris@math.uc3m.es}

\maketitle{}


\begin{abstract}
In this paper we prove that any Riemannian surface,
with no restriction of curvature at all,
can be decomposed into blocks belonging just to some of these types:
generalized Y-pieces, generalized funnels and halfplanes.\\
\hspace{-.6cm} \textit{Key words and phrases:} Decomposition of surfaces; arbitrary
curvature.\\

\end{abstract}

\section{Introduction.}

The Classification Theorem of compact surfaces
says that every orientable compact
topological surface is homeomorphic either to a sphere
or to a ``torus" of genus $g \ge 1$
(see e.g. \cite{M}).

\smallskip

We say that the closure of a three-holed sphere (which is a bordered compact
topological surface whose border is the union of three pairwise disjoint
simple closed curves) is a \emph{topological Y-piece}.
A Y-piece can be visualized as a tubing with the shape of the letter Y.
A \emph{cylinder} is a bordered topological surface
homeomorphic to $\SS^1\times [0,\infty)$,
where $\SS^1$ is the unit circle.

We refer to the next section for precise definitions and background.

\smallskip

The Classification Theorem of compact surfaces
says, in other words, that every orientable compact
topological surface except for the sphere and the torus (of genus $1$)
can be obtained by gluing topological Y-pieces along their boundaries.

\smallskip

In \cite{AR}, the Classification Theorem is generalized to noncompact
surfaces in the following way:

\smallskip

\begin{teorema}
\label{t:descomposicion02}
$($\cite[Theorem 1.1]{AR}$)$
Every complete orientable topological surface
which is homeomorphic neither to the sphere nor to the plane nor to the torus
is the union (with pairwise disjoint interiors)
of topological Y-pieces and cylinders.
\end{teorema}

\smallskip

The following result is the most important
in \cite{AR} and is a geometric version of the theorem above for complete
surfaces with constant negative curvature.
In this case we have more information about the basic blocks of the
surface: the surface can be decomposed in such a way that the boundary
of the blocks is the union of at most three simple closed geodesics.
Since the Riemannian structure is more restrictive than the
topological one, an additional piece is necessary in order to achieve the
decomposition: the halfplane.

We state now this result for Riemannian surfaces.

\smallskip

\begin{teorema}
\label{t:descomposicion01}
$($\cite[Theorem 1.2]{AR}$)$
Every complete orientable Riemannian surface with constant curvature $K=-1$,
which is not the punctured disk,
is the union (with pairwise disjoint interiors)
of generalized Y-pieces, funnels and halfplanes.
\end{teorema}

\smallskip

In the applications of Theorem \ref{t:descomposicion01},
it is a crucial fact that the
boundaries of the generalized Y-pieces are simple closed geodesics.
There is a clear reason for this: it is very easy to cut and paste
surfaces along such kind of curves.

Furthermore, closed geodesics in a Riemannian surface $S$
are geometrical objects interesting by themselves.
Since they are the periodic orbits of the dynamical system associated to
$S$ on its unit tangent bundle, they provide tools to study the geodesic
flow, just like the fixed points of an automorphism helps to study it.
Lastly, the closed geodesics are becoming more and more important in the study
of heat and wave equations, and of the spectrum of $S$.
The lengths of all closed geodesics determine largely the spectrum.
Conversely, the spectrum determines completely
the lengths of the closed geodesics
(see \cite{Ch}, \cite{G}, \cite{CV}).

\smallskip

In this paper we prove the conclusion of Theorem
\ref{t:descomposicion01} with no restriction of curvature at all
(see Theorem \ref{t:descomposicion} forward). In our context, we
require that the boundaries of the generalized Y-pieces are
minimizing simple closed geodesics (in its free homotopy class).
Although if $K=-1$ the property of minimization ever holds (see e.g.
Theorem \ref{t:geodconknegativa} and Lemma \ref{l:unicidad}),
this is obviously false for arbitrary curvature.

Theorem \ref{t:descomposicion2} is a sharp version of
Theorem \ref{t:descomposicion}
for Riemannian surfaces with curvature $K\le -c^2 < 0$.

Finally,
Theorems \ref{t:descomposicion3} and \ref{t:descomposicion4}
are versions, respectively, of
Theorems \ref{t:descomposicion} and \ref{t:descomposicion2}
in the context of bordered Riemannian surfaces.

\smallskip

J. L. Fern\'andez and M. V. Meli\'an (see \cite{FM}) proved the
following result which helps to understand the behavior of geodesics
in surfaces with curvature $K=-1$.

\smallskip

\begin{teorema}
\label{t:FM}
$($\cite[Theorem 1]{FM}$)$
Let $S$ be a complete orientable Riemannian surface with curvature $K=-1$.
There are three possibilities:

$(i)$ $S$  has finite area. Then for every $p\in S$ there is
exactly a countable collection of directions in $\Ep(p)$.

$(ii)$ $S$ is transient. Then for every $p\in S$, $\Ep(p)$ has full
measure.

$(iii)$ $S$ is recurrent and of infinite area.
Then for every $p\in S$, $\Ep(p)$ has zero
length but its Hausdorff dimension is $1$.
\end{teorema}

\smallskip

A surface is said to be \emph{transient}
(respectively, \emph{recurrent})
if Brownian motion of $S$ is transient (respectively, recurrent).
Also, we define $\Ep(p)$ as the set of unitary directions $v$ in the tangent
plane of $S$ at $p$ such that the unit speed geodesic emanating from
$p$ in the direction of $v$, \emph{escapes} to infinity.

\smallskip

Just like
Theorem \ref{t:descomposicion01}
played an important role in the proof of Theorem \ref{t:FM},
we are sure that
Theorems \ref{t:descomposicion} and
\ref{t:descomposicion2} will be crucial in order to generalize
this latest result to surfaces with curvature $K \le -c^2 < 0$.

\smallskip

The argument in the proof of
Theorem \ref{t:descomposicion} is quite alike to the one in the proof of
Theorem \ref{t:descomposicion01}.
Unfortunately, \emph{every} standard fact used in the proof of
Theorem \ref{t:descomposicion01} is false when there is no
restriction of curvature.
Hence, it was unavoidable
both to state definitions for the new objects appearing in our current context
and to prove alternate results valid for arbitrary curvature.
This work has provided some results with intrinsic interest, as
Theorems \ref{t:geodconknegativa} and \ref{t:disjuntas}.

\smallskip

One can think that in the decomposition of
Theorem \ref{t:descomposicion}
we might not need halfplanes.
There is an example in \cite{AR} which shows that,
even with curvature $K=-1$, we do need them.
The necessity of halfplanes is, in fact, one of the most
difficult parts in the proof of this theorem.

\smallskip

The outline of the paper is as follows.
Section 2 presents some definitions and technical results which we will need.
We prove some additional technical results in Section 3.
Section 4 is dedicated to the main results.

\smallskip

\noindent {\bf Notations.} We denote by $L_M(\g)$ the length of a
curve $\g$ in a Riemannian manifold $M$. If there is no possible
confusion, we usually write $L(\g)$.

\smallskip

\noindent {\bf Acknowledgements.} We would like to thank
Professor Jes\'us Gonzalo his proof of
Theorem \ref{t:disjuntas}.

\section{Background in Riemannian manifolds.}

\begin{definicion}
Any divergent curve $\s: [0,\infty) \longrightarrow Y$,
where $Y$ is a noncompact Hausdorff space, determines an
\emph{end} $E$ of $Y$.
Given a compact set $F$ of $Y$,
one defines $E(F)$ to be the arc component of $Y \setminus F$
that contains a terminal segment $\s([a,\infty))$ of $\s$.
A set $U \subset Y$ is a
\emph{neighborhood of an end} $E$ if
$U$ contains $E(F)$ for some compact set $F$ of $Y$.
An end $E$ in a surface $S$ is \emph{collared}
if $E$ has a neighborhood homeomorphic to
$(0,\infty) \times \SS^1$.
A neighborhood $U$ of $E$ will be called
\emph{Riemannian collared}
if there exists a $C^1$ diffeomorphism
$X: (0,\infty) \times \SS^1 \longrightarrow U$
such that the metric in $U$ relative to the coordinate system
$X$ is written $ds^2= dr^2 + G(r,\th)^2 d\th^2$,
where $G$ is a positive continuous function.
A sequence of curves $\{C_n\}$
\emph{converges to} $E$
if for any neighborhood $U$ of $E$
we have $C_n \subset U$ for sufficiently large $n$.
We say that a closed curve $\g$ \emph{bounds} a collared end $E$
in $S$ if some arc component of $S \setminus \g$
is a neighborhood $U$ of $E$.
\end{definicion}

It follows directly from the metric expression
$ds^2= dr^2 + G(r,\th)^2 d\th^2$
of a Riemannian collared parametrization
that the $r$-parameter curves have unit speed and
minimize the distance between any two of their points.
Consequently the $r$-parameter curves are geodesics of $S$.
If the curvature $K$ satisfies
$K \le 0$, then $G$ is a $C^{\infty}$ function of $r$
for each fixed $\th$
and satisfies the Jacobi equation
$$
\frac{\p^2 G}{\p r^2}\, (r,\th)
+ K(r,\th) G(r,\th)=0\,,
$$
where $K(r,\th)$ is the curvature of $S$ at $X(r,\th)$
(\cite[p. 17]{E}).

\medskip

Every manifold is connected, $C^{\infty}$ and satisfy the second
axiom of countability (has a countable basis for its topology). In a
Riemannian surface we always assume that the Riemannian metric is
$C^{\infty}$ unless perhaps in some simple closed geodesics, each of
them bounding a collared end, where we allow the metric to be
$C^{1}$ and piecewise $C^{\infty}$, with the ``singularities" along
these geodesics. Then the curvature is a (possibly discontinuous)
function along these geodesics.

Geodesic always means local geodesic (unless we say explicitly
something else).

\begin{definicion}
Given a Riemannian surface $S$, a geodesic $\g$ in $S$,
and a continuous unit vector field $\xi$ along $\g$,
orthogonal to $\g$, we define the
\emph{Fermi coordinates} based on $\g$ as the map
$Y(r,\th):=\exp_{\g(\th)} r \xi(\th)$.
\end{definicion}

It is well known that the Riemannian metric can be expressed in
Fermi coordinates as $ds^2= dr^2 + G(r,\th)^2 \, d\th^2$, where
$G(r,\th)$ is the solution of the scalar equation
$$
\frac{\p^2 G}{\p r^2}\, (r,\th)
+ K(r,\th) G(r,\th)=0\,,
\qquad G(0,\th)=1\,,
\qquad \frac{\p G}{\p r}\, (0,\th)=0\,,
$$
(see e.g. \cite[p. 247]{C}).

We will need the following three results.

\begin{teorema}
\label{t:ends1}
$($\cite[Theorem 4.2]{E}$)$
Let $S$ be a complete Riemannian surface with $K \le 0$
and $E$ an end of $S$.
Then the following are equivalent:

$(1)$ $E$ is a collared end.

$(2)$ $E$ is a Riemannian collared end.

$(3)$ There exists a sequence $\{C_n\}$ of continuous piecewise smooth
closed curves converging to $E$ such that $\{C_n\}$ belongs to
a single nontrivial free homotopy class.
\end{teorema}

It is clear that $(2)$ can be deduced from $(1)$
since $S$ verifies $K \le 0$.
However, since $(1)$ and $(3)$ are topological conditions,
we have the following result without conditions on $K$.

\begin{teorema}
\label{t:ends2}
Let $S$ be a complete Riemannian surface
and $E$ an end of $S$.
Then $E$ is a collared end if and only if
there exists a sequence $\{C_n\}$ of continuous piecewise smooth
closed curves converging to $E$ such that $\{C_n\}$ belongs to
a single nontrivial free homotopy class.
\end{teorema}

\begin{teorema}
\label{t:metricspaces}
$($\cite[Theorem (5.16)]{B}$)$
Giving a sequence of rectifiable curves $\{\a_k\}$
contained in a compact set of a Riemannian manifold $M$
with $\{L(\a_k)\}$ a bounded sequence,
there exists a subsequence of curves
(which we also call $\{\a_k\}$ for simplicity),
a rectifiable curve $\a$,
and parametrizations $x_k:[0,1]\longrightarrow X$ of $\a_k$
and $x:[0,1]\longrightarrow X$ of $\a$,
such that $\{x_k\}$ converges uniformly to $x$ in $[0,1]$
and
$$
L(\a) \le \liminf_{k \to \infty} L(\a_k) \,.
$$
\end{teorema}

In fact, Theorem (5.16) in \cite{B} is stronger than
Theorem \ref{t:metricspaces}, but this statement
is general enough for our purposes.

\begin{definicion}
Given a $n$-dimensional Riemannian manifold $M$
and a closed curve $\a$ in $M$, we define the
\emph{length of the freely homotopy class} $[\a]$ as
$$
L([\a]) := \inf \big\{ L(\s) : \, \s \in [\a] \big\} \,.
$$
\indent
The curve $\a$ is \emph{minimizing} if $L(\a) = L([\a])$.

A \emph{minimizing sequence for $\a$}
is a sequence of closed curves $\{\a_k\} \subset [\a]$
such that $\lim_{k \to \infty} L(\a_k) = L([\a])$.
\end{definicion}

\begin{definicion}
A \emph{halfplane} is a bordered Riemannian surface
which is simply connected
and whose border is a unique nonclosed simple geodesic.

A \emph{generalized funnel} is a bordered Riemannian surface
which is a neighborhood of a collared end
and whose border is a minimizing simple closed geodesic.
A \emph{funnel} is a generalized funnel
such that there does not exist another simple closed geodesic
freely homotopic to the border of the funnel.

\begin{figure}[h]
    \centering
    \subfigure[Generalized funnel]{\epsfig{figure=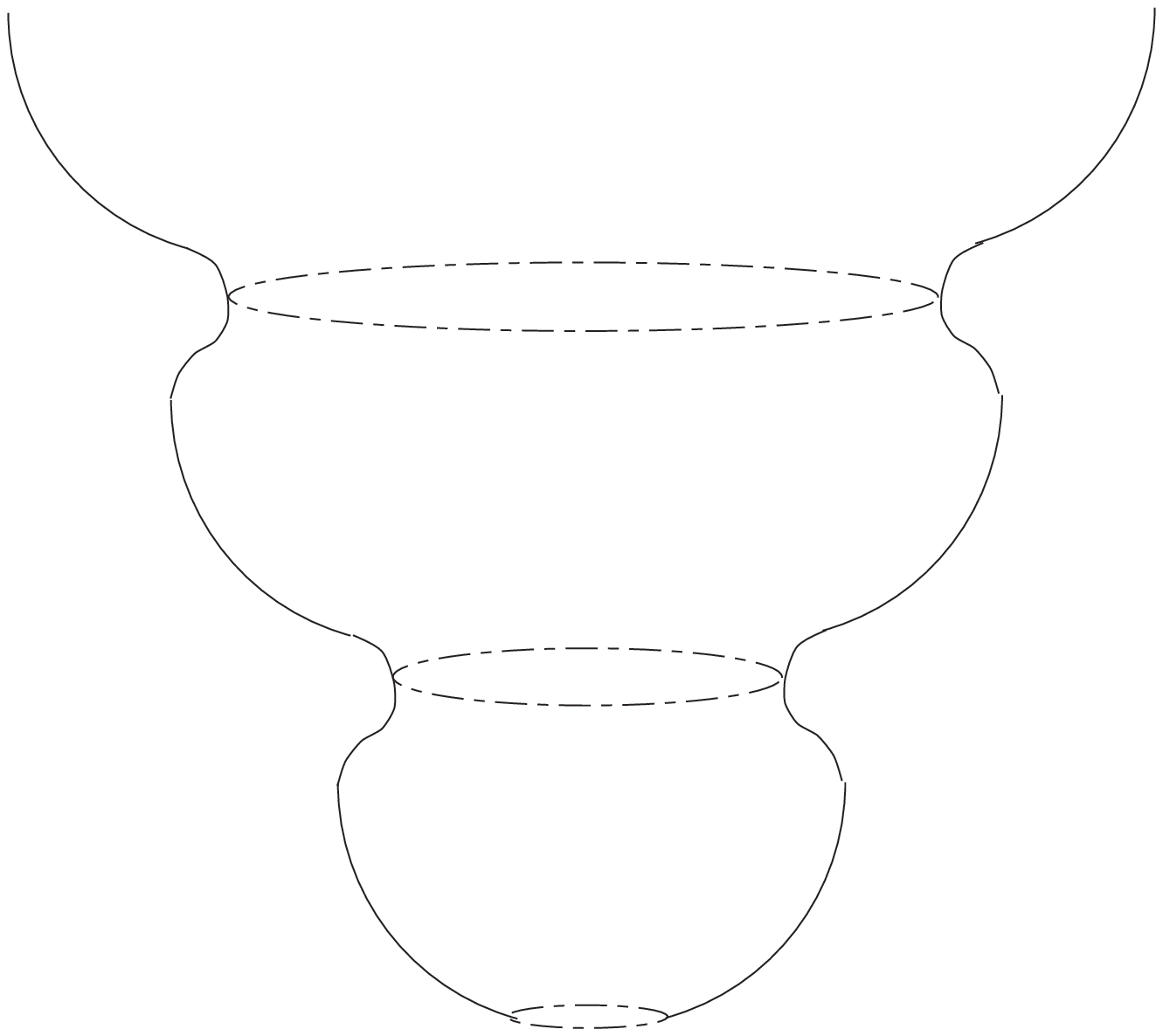,height=4cm}}\qquad\qquad %
    \subfigure[Generalized funnel]{\epsfig{figure=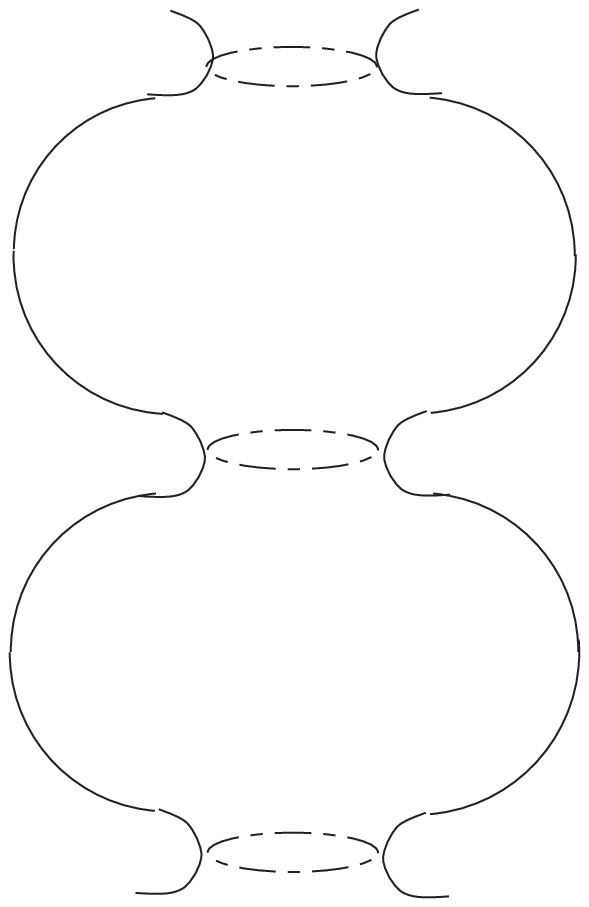,width=4cm,height=4cm}}\qquad\qquad %
    \subfigure[Funnel]{\epsfig{figure=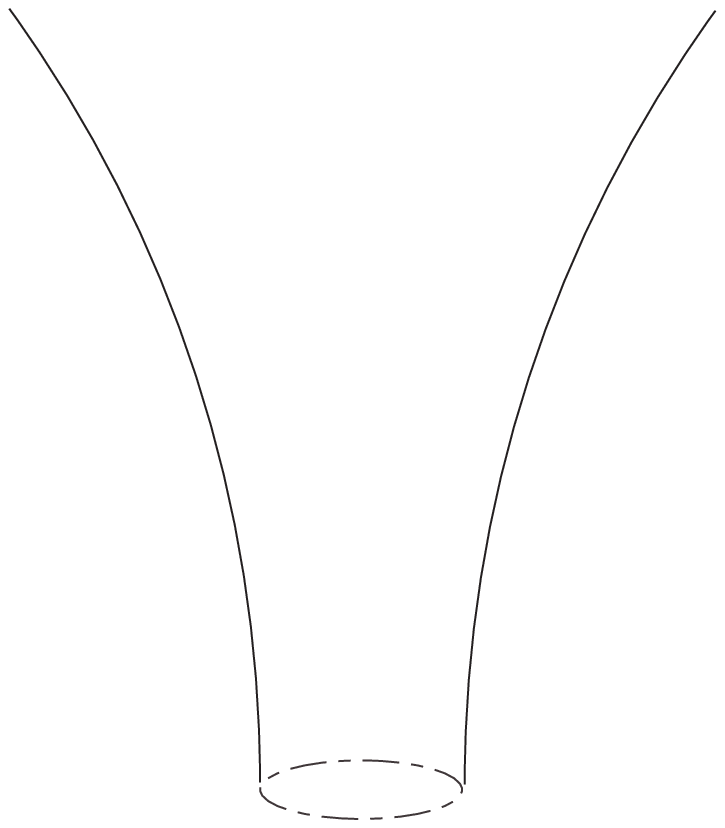,height=4cm}}
\end{figure}

A \emph{generalized puncture} is a collared end whose fundamental
group is generated by a simple closed curve $\s$ and there is no
minimizing closed geodesic $\g \in [\s]$. A \emph{puncture} is a
generalized puncture such that $L([\s])=0$ and there is no closed
geodesic in $[\s]$.

\begin{figure}[h]
    \centering
    \subfigure[$\!\!$Generalized puncture]{\epsfig{figure=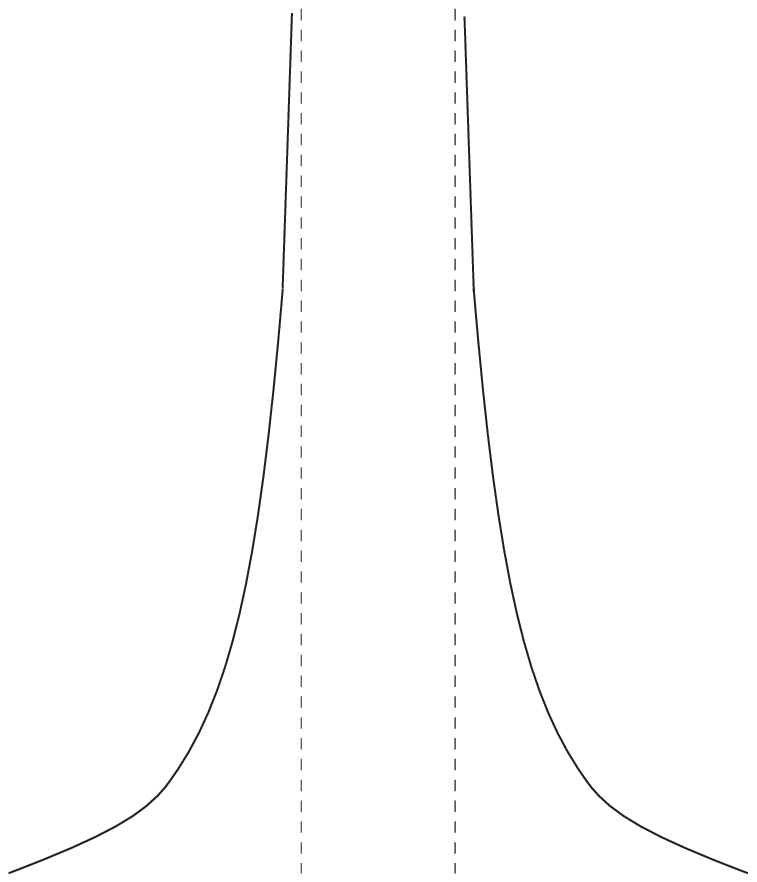,height=4.5cm}}\qquad\qquad %
    \subfigure[Generalized puncture]{\epsfig{figure=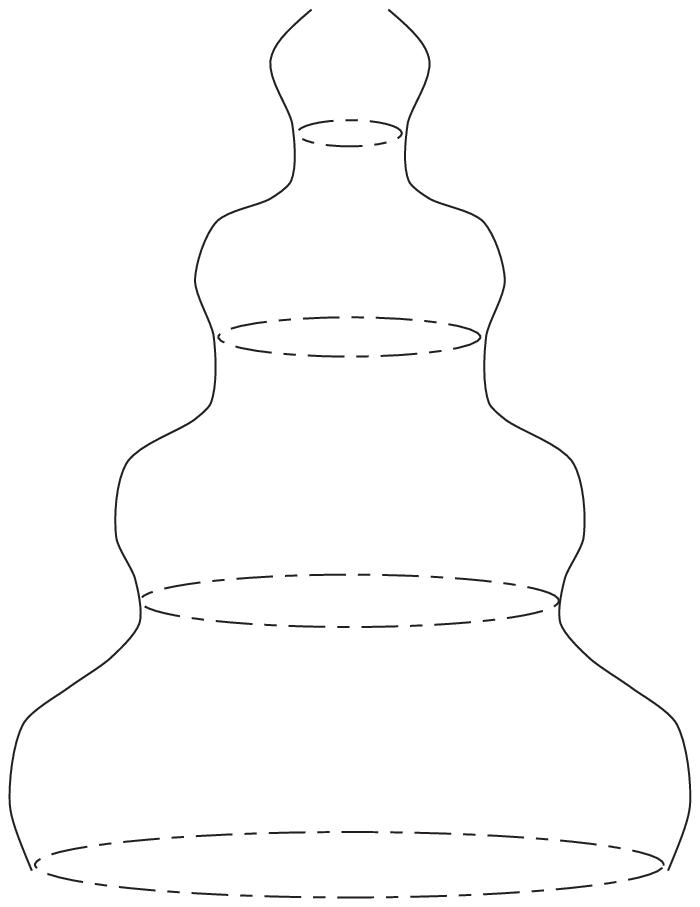,width=4cm,height=4cm}}\qquad\qquad %
    \subfigure[Puncture]{\epsfig{figure=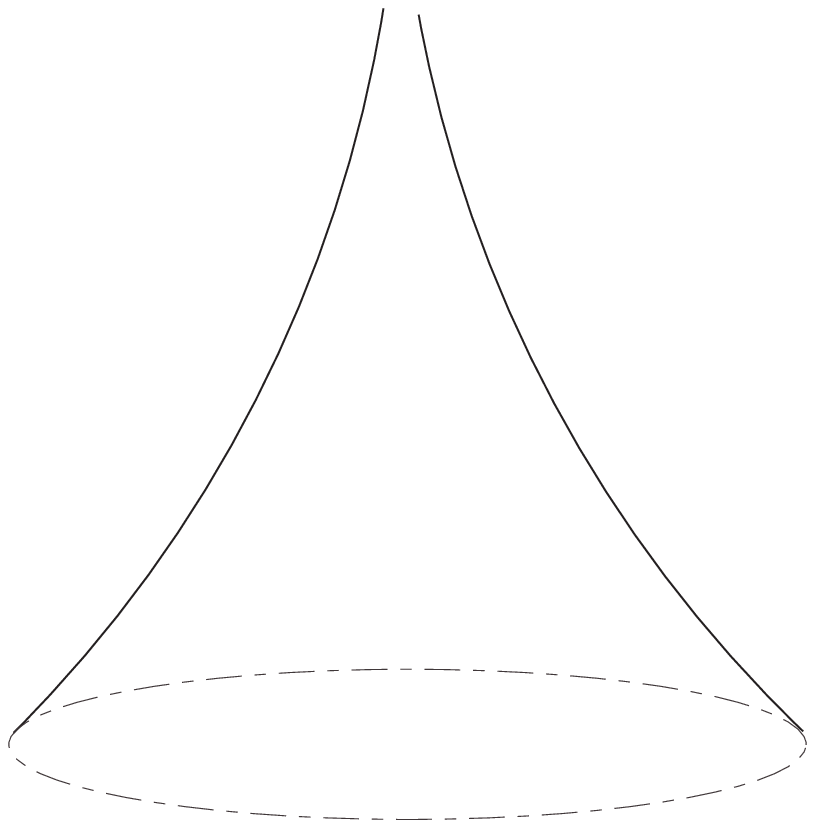,height=4cm}}%
\end{figure}

A bordered or nonbordered surface is \emph{doubly connected}
if its fundamental group is isomorphic to $\ZZ$.
Every generalized funnel and every generalized puncture are
doubly connected surfaces.

A \emph{geodesic domain} $G$ is a bordered Riemannian surface
(which is neither simply nor doubly connected)
with finitely generated fundamental group and
such that $\p G$ consists of finitely many minimizing simple closed geodesics,
and it may contain generalized punctures but not generalized funnels.

A \emph{Y-piece} is a compact bordered Riemannian
surface which is topologically a sphere without three open disks
and whose boundary curves are minimizing simple closed geodesics.
They are a standard tool to construct Riemannian surfaces.
A clear description of these Y-pieces and their use is given
in \cite[chapter X.3]{C}.

A \emph{generalized Y-piece} is a bordered or nonbordered
Riemannian surface which is topologically a sphere without three open
disks, such that there exist integers $n,m\ge 0$ with $n+m=3$,
so that the border are $n$ minimizing
simple closed geodesics and there are $m$ generalized punctures.


Notice that a generalized Y-piece is topologically
the union of a Y-piece and $m$ cylinders, with $0\le m\le 3$.
It is clear that every generalized Y-piece
is a geodesic domain (unless $m=3$, in which case it has no border).
Furthermore, every geodesic domain is a finite
union (with pairwise disjoint interiors) of generalized Y-pieces
(see Proposition \ref{p:domgeod} below).

We say that the set $A$ is \emph{exhausted} by $\{A_n\}$ if $A_n\subseteq
A_{n+1}$ for every $n$ and $A=\cup_n A_n$.

We say that a bordered Riemannian surface $S$ is \emph{simple}
if the border of $S$ is a (finite or infinite)
union of pairwise disjoint simple closed geodesics.
\end{definicion}

\section{Technical results.}

\begin{lema}
\label{l:exists}
Let us consider a $n$-dimensional complete Riemannian manifold $M$
and an homotopically nontrivial closed curve $\a$ in $M$.
If there exists a minimizing sequence $\{\a_k\}$
for $\a$ contained in a compact set,
then there exists a minimizing closed geodesic $\g \in [\a]$.
\end{lema}

\begin{proof}
Since $\{L(\a_k)\}$ is convergent,
it is a bounded sequence.
By Theorem \ref{t:metricspaces},
there exists a subsequence of curves
(which we also call $\{\a_k\}$ for simplicity),
a rectifiable curve $\g$,
and parametrizations $x_k:[0,1]\longrightarrow M$ of $\a_k$
and $x:[0,1]\longrightarrow M$ of $\g$,
such that $\{x_k\}$ converges uniformly to $x$ in $[0,1]$
and
$$
L(\g) \le \liminf_{k \to \infty} L(\a_k) = L([\a]) \,.
$$
The curve $\g$ is closed since each $\a_k$ is a closed curve and
$x(0) = \lim_{k\to\infty} x_k(0) = \lim_{k\to\infty} x_k(1) = x(1)$.
Then, in order to finish the proof of the lemma, it is enough to show that
$\g \in [\a]$,
since then $\g$ attains the minimum length in its homotopy class,
and it must be a geodesic.

We can assume that $x_k$ and $x$ are $1$-periodic functions in $\RR$.
For each $t \in [0,1]$ let us consider $r_t>0$ small enough
to guarantee that the ball $B(x(t),r_t)$ in $M$
is simply connected.
For each $t \in [0,1]$, let us denote by $J_t$
the connected component of
$\g \cap B(x(t),r_{t})$ which contains $x(t)$.
Since $\g$ is a compact topological space and
$\{J_t\}_{t\in [0,1]}$ is an open covering of $\g$,
there exist
$0\le t_1< t_2 < \cdots < t_{m-1} < t_m \le 1$
such that
$\g \subset \cup_{j=1}^m J_{t_j}$.
Choosing a subset of
$\{t_1, t_2, \dots , t_m \}$
if it is necessary,
without loss of generality we can assume that the subcovering
$\{J_{t_j}\}_{j=1}^m$
is minimal in the following sense:
each $y\in \g$ belongs at most to two sets of $\{J_{t_j}\}_{j=1}^m$.

Let us consider
$0\le s_1< \cdots < s_{m-1} < 1$  and $s_m \in (s_{m-1}, s_1 + 1)$
such that
$$
x(s_1)\in J_{t_1} \cap J_{t_2}, \dots,
x(s_{m-1})\in J_{t_{m-1}} \cap J_{t_m},
x(s_{m})\in J_{t_m} \cap J_{t_1} \,.
$$
Hence,
$$
x(s_{m}), x(s_1) \in B(x(t_1),r_{t_1}), \;
x(s_{1}), x(s_2) \in B(x(t_2),r_{t_2}), \; \dots \, , \;
x(s_{m-1}), x(s_m) \in B(x(t_m),r_{t_m}).
$$
Since $\{x_k\}$ converges uniformly to $x$ in $\RR$,
there exists $k_0$ such that
$$
x_k([s_m,s_1])\subset B(x(t_1),r_{t_1}), \;
x_k([s_1,s_2])\subset B(x(t_2),r_{t_2}), \; \dots \, , \;
x_k([s_{m-1},s_m])\subset B(x(t_m),r_{t_m}),
$$
for every $k\ge k_0$.

This proves that $\g \in [\a_k]$ for every $k\ge k_0$
(since each ball $B(x(t),r_t)$ is simply connected),
and then $\g \in [\a]$.
This finishes the proof of the lemma.
\end{proof}

\begin{prop}
\label{p:alternativa}
Let us consider a complete Riemannian surface $S$
and an homotopically nontrivial closed curve $\a$ in $S$.
Then, one and only one of the two following possibilities holds:

$(1)$ There exists a minimizing closed geodesic $\g \in [\a]$.

$(2)$ The curve $\a$ bounds a generalized puncture $E$ in $S$.
Furthermore, if $S$ is not doubly connected,
then any minimizing sequence for $\a$ converges to $E$.
\end{prop}

\begin{proof}
If there exists a minimizing sequence $\{\a_k\}$
for $\a$ contained in a compact set,
then Lemma \ref{l:exists} gives $(1)$.

Otherwise, every minimizing sequence $\{\a_k\}$ for $\a$
escapes from any compact set.
Hence, there not exists any minimizing closed geodesic in $[\a]$.

If $S$ is doubly connected, then $\a$ bounds a collared end $E$
(in fact, $\a$ bounds exactly two collared ends).
This collared end $E$ is a generalized puncture since
there not exists any minimizing closed geodesic in $[\a]$.

If $S$ is not doubly connected, then
any minimizing sequence $\{\a_k\}$ for $\a$
converges to an end $E$.
Since the curves $\{\a_k\}$ belong to
a single nontrivial free homotopy class,
Theorem \ref{t:ends2} gives that $E$
is a collared end in $S$.
Hence, $\a$ bounds a collared end in $S$,
which must be a generalized puncture since
there not exists any minimizing closed geodesic in $[\a]$.
\end{proof}

In order to deal with bordered surfaces, we need the following results.

\begin{lema}
\label{l:complecion0}
Any simple complete bordered Riemannian surface
$S$ is a subset of a complete Riemannian surface $R$, which can be
obtained by attaching a neighborhood of a collared end
to each simple closed geodesic $\g\subseteq\partial S$, with the following properties:

$(1)$ If $\s$ is a closed curve in $R$ which is not contained in $S$, then there exists
$\s_0\subseteq (S\cap \s)\cup \p S \subset S$ with $\s_0\in [\s]$ and $L(\s_0) < L(\s)$.

$(2)$ A closed geodesic is minimizing in $R$ if and only if it is minimizing in $S$
(in particular, it is contained in $S$).

$(3)$ If $\s$ is a closed curve in $S$, then $L_S([\s])=L_R([\s])$.

$(4)$ If $\s$ is a closed curve in $S$, and
$\{\s_k\}$ is a minimizing sequence for $\s$ verifying
$\{\s_k\} \subseteq (R \setminus S)\cup K$, with $K$ a compact subset of $S$,
then there exists a minimizing closed geodesic in $[\s]$.

$(5)$ The curvature satisfies $K=-1$ in $R \setminus S$.

$(6)$ The fundamental group of $R$ is isomorphic to the fundamental group of $S$.

$(7)$ If $S$ is not doubly connected, then
there exists a minimizing simple closed geodesic in $[\g_0]$
for each simple closed geodesic $\g_0 \subseteq \p S$.
\end{lema}

\begin{proof}
The border of $S$ is a (finite or infinite)
union of pairwise disjoint simple closed geodesics.
Let us fix a closed geodesic $\g_0 \subseteq \p S$ with length $l$.
We can consider the Fermi coordinates based on $\g_0$.
The Riemannian metric can be expressed in
Fermi coordinates as $ds^2= dr^2 + G(r,\th)^2 \, d\th^2$, with
$G(r,\th)$ a $l$-periodic function in $\th$
defined in $[-r_0,0] \times \RR$, for some $r_0 > 0$.
We have $G(0,\th)=1$ and $\p G/\p r (0,\th)=0$ for every $\th \in \RR$.
If we define $G(r,\th):= \cosh r$
in $(0,\infty) \times \RR$, then
it is $C^1$ (and even piecewise $C^{\infty}$) in $[-r_0,\infty) \times \RR$,
and $l$-periodic in $\th$.
These coordinates
$(r,\th) \in [-r_0,\infty) \times \RR$,
with the Riemannian metric
$ds^2= dr^2 + G(r,\th)^2 \, d\th^2$,
attach a neighborhood of a collared end $F$ to $\g_0$;
by this way we get a $C^{\infty}$ surface.
We have that $K(r,\th)=-1$ in $(0,\infty) \times \RR$.
We also have the following properties:

$(a)$ Any homotopically nontrivial closed curve $\s$ in $F$ verifies
$L(\g_0) < L(\s)$:

Without loss of generality we can assume that
$\s$ can be parametrized in Fermi coordinates based on $\g_0$ as
$\s(\th)=(r(\th),\th)$, with $\th \in [0,l]$.
Then,
$$
L(\s) = \int_0^l \sqrt{r'(\th)^2 + \cosh\!^2 r(\th)} \; d\th
\ge \int_0^l \cosh r(\th)  \, d\th
> \int_0^l d\th = l = L(\g_0) \,.
$$

$(b)$ Given any closed curve $\s$ intersecting $S$ and the interior of $F$,
there exists $\s_0 \in [\s]$ contained in $S$ verifying
$L(\s_0) < L(\s)$:

We can construct this curve in the following way:
given any subarc $a$ of $\s$ contained in $F$ and joining two points $p,q \in \g_0$,
we replace it by the subarc of $\g_0$ joining $p,q,$
which is homotopic to $a$.
The argument above gives $L(\s_0) < L(\s)$.

We define $R$ as the surface obtained by
attaching this neighborhood of a collared end to each
closed geodesic in $\p S$.

Properties $(a)$ and $(b)$ give that
if $\s$ is a closed curve in $R$ which is not contained in $S$, then there exists
$\s_0\subseteq (S\cap \s)\cup \p S \subset S$
with $\s_0\in [\s]$ and $L(\s_0) < L(\s)$.
This finishes the proof of $(1)$, $(5)$ and $(6)$.

Now, the statements $(2)$ and $(3)$ are direct consequences of $(1)$.

We prove now $(4)$.
If $\s$ is a closed curve in $S$, and
$\{\s_k\}$ is a minimizing sequence for $\s$ verifying
$\{\s_k\} \subseteq (R \setminus S)\cup K$, with $K$ a compact subset of $S$,
by $(1)$ there exists $\{\s_k^0\} \subseteq K$
with $\s_k^0\in [\s]$ and $L(\s_k^0) \le L(\s_k)$.
Then $\{\s_k^0\}$ is a minimizing sequence for $\s$ contained in a compact set
and Lemma \ref{l:exists} gives
that there exists a minimizing closed geodesic in $[\s]$.

In order to prove $(7)$, fix a simple closed geodesic $\g_0
\subseteq \p S$. Let us call $F$ the neighborhood of a collared end
in $R$ with $\p F=\g_0$. Seeking for a contradiction, assume that
there not exists a minimizing closed geodesic in $[\g_0]$. Since $R$
is not doubly connected, by Proposition \ref{p:alternativa} $\g_0$
bounds a generalized puncture $E$ in $R$ and any minimizing sequence
for $\a$ converges to $E$. Since $R$ is not doubly connected, $F$ is
a neighborhood of $E$, and for any minimizing sequence $\{\a_k\}$
for $[\g_0]$ there exists $N$ with $\a_k\subset F$ for every $k\ge
N$. By $(a)$ we have $L(\g_0) < L(\a_k)$ for every $k\ge N$, which
is the required contradiction.
\end{proof}

Using Lemma \ref{l:complecion0},
Proposition \ref{p:alternativa}
can be generalized to simple bordered Riemannian surfaces.

\begin{prop}
\label{p:alternativaborde}
Let us consider a simple complete bordered Riemannian surface $S$
and an homotopically nontrivial closed curve $\a$ in $S$.
Then, one and only one of the two following possibilities holds:

$(1)$ There exists a minimizing closed geodesic $\g \in [\a]$.

$(2)$ The curve $\a$ bounds a generalized puncture $E$ in $S$
and any minimizing sequence for $\a$ converges to $E$.
\end{prop}

\begin{proof}
By Lemma \ref{l:complecion0},
$S$ is a subset of a complete Riemannian surface $R$.

Assume first that $S$ is not doubly connected
(then $R$ is not doubly connected).

By Lemma \ref{l:complecion0} (2),
if there exists a minimizing closed geodesic $\g \in [\a]$ in $R$,
then $\g \in S$.

If there not exist such minimizing geodesic,
by Proposition \ref{p:alternativa}
the curve $\a$ bounds a generalized puncture $E$ in $R$
and any minimizing sequence for $\a$ converges to $E$.
By Lemma \ref{l:complecion0} (7),
$\a$ can not be freely homotopic to any closed geodesic in $\p S$,
and therefore $E \subset S$.

Assume now that $S$ is doubly connected
(then $R$ is also doubly connected).
The curve $\a$ bounds exactly two collared ends in $R$.

Since $S$ is doubly connected, $\p S$ can be either
a simple closed geodesic or
two simple closed geodesics.

Assume first that $\p S$ is a simple closed geodesic $\g_0$.
Then $\a$ bounds the unique collared end $E$ in $S$.
Consider a minimizing sequence $\{\a_n\}$ for $[\a]$ in $R$.
If $(2)$ does not hold, then
either $\a$ does not bound a generalized puncture
(and then $(1)$ holds) or
there exists a neighborhood $U$ of $E$
and a subsequence $\{\a_{n_k}\}$ with
$\a_{n_k}\nsubseteq U$ for every $k$.
Since $\{L(\a_{n})\}$ is a bounded sequence,
without loss of generality we can assume that
$\a_{n_k}\cap U= \varnothing$ for every $k$.
Since $S$ is doubly connected, $R=S\cup F$ with $\p F= \p S=\g_0$, and
the complement of $U$ in $R$ is $F\cup K$, where $K$ is a compact set in $S$.
Therefore,
$\a_{n_k}\subset F\cup K$ for every $k$ and by Lemma \ref{l:complecion0} $(4)$
there exists a minimizing closed geodesic in $[\a]$.
Then $(1)$ also holds.

Assume now that $\p S$ is the union of two simple closed geodesics $\g_1,\g_2$.
Then, $R=S\cup F_1 \cup F_2$ with $\p F_j= \g_j$ ($j=1,2$)
and $S$ is compact.
Consider a minimizing sequence $\{\a_n\}$ for $[\a]$ in $R$.
By Lemma \ref{l:complecion0} $(1)$, there exists
a minimizing sequence $\{\a_n^0\}\subset S$ for $[\a]$.
Since $S$ is compact,
Lemma \ref{l:exists} gives that there exists
a minimizing closed geodesic $\g \in [\a]$.
\end{proof}

\begin{lema}
\label{l:sturm}
Let us consider a positive constant $c$ and two functions
$y_0,y,$ satisfying respectively
$y_0''=c^2 y_0$, $y_0(t_0) > 0$, $y_0'(t_0) > 0$,
and
$$
\left\{ \begin{array}{ll}
y''(t)\ge c^2 y(t) >0 \,,
& \;  \text{ if } \, t\ge t_0 \,,
\\
y(t_0)= y_0(t_0) \,,
& \;
\\
y'(t_0) \ge y_0'(t_0) \,.
& \;
\end{array}
\right.
$$
Then, $y(t)\ge y_0(t)$ for every $t\ge t_0$.
\end{lema}

\begin{proof}
Since $y''(t)>0$ if $t\ge t_0$, then $y'$ is an increasing function in $[t_0,\infty)$.
This fact and $y'(t_0) \ge y_0'(t_0) > 0$ give
$y'(t)\ge y'(t_0) > 0$ for every $t\ge t_0$.
Then, for every $t\ge t_0$, we can deduce
$$
\aligned
y''(t) & \ge c^2 y(t) \,,
\\
y''(t) y'(t) & \ge c^2 y(t) y'(t) \,,
\\
y'(t)^2 - y'(t_0)^2 & \ge c^2 \big( y(t)^2 - y(t_0)^2 \big) \,,
\\
y'(t) & \ge \sqrt{ c^2 \big( y(t)^2 - y(t_0)^2 \big) + y'(t_0)^2} \;.
\endaligned
$$
For each fixed $\e \in (0,y_0'(t_0))$, we define the function
$y_{\e}$ as the unique solution of
$$
\left\{ \begin{array}{ll}
y_{\e}''(t)= c^2 y_{\e}(t)  \,,
& \;  \text{ if } \, t\ge t_0 \,,
\\
y_{\e}(t_0)= y_0(t_0) \,,
& \;
\\
y_{\e}'(t_0) = y_0'(t_0)-\e > 0 \,.
& \;
\end{array}
\right.
$$
Then $y'(t_0)> y_{\e}'(t_0) > 0$.
Using the same argument above in the case of $y_{\e}$ (with equality instead of inequality)
we obtain that
$y_{\e}'(t)\ge y_{\e}'(t_0) > 0$ for every $t\ge t_0$ and
$$
y_{\e}'(t) = \sqrt{ c^2 \big( y_{\e}(t)^2 - y_{\e}(t_0)^2 \big) + y_{\e}'(t_0)^2} \;.
$$
We prove now that
$y(t)\ge y_{\e}(t)$ for every $t\ge t_0$.
Seeking for a contradiction, suppose that
$y(t)< y_{\e}(t)$ for some $t> t_0$.
Then, we can define $t_1:=\min\{ t> t_0 : \, y(t) = y_{\e}(t) \}$;
this minimum is attained since $y(t_0)= y_{\e}(t_0)$ and $y'(t_0)> y_{\e}'(t_0)$;
consequently, $y(t) > y_{\e}(t) > 0$ for every $t\in (t_0,t_1)$, and
$$
\aligned
& y_{\e}(t_1) - y_{\e}(t_0)
= \int_{t_0}^{t_1} y_{\e}'(t) \, dt
= \int_{t_0}^{t_1} \sqrt{ c^2 \big( y_{\e}(t)^2 - y_{\e}(t_0)^2 \big) + y_{\e}'(t_0)^2} \, dt
\\
& < \int_{t_0}^{t_1} \sqrt{ c^2 \big( y(t)^2 - y(t_0)^2 \big) + y'(t_0)^2} \, dt
\le \int_{t_0}^{t_1} y'(t) \, dt
= y(t_1) - y(t_0)
= y_{\e}(t_1) - y_{\e}(t_0) \,.
\endaligned
$$
This is a contradiction and we have proved that
$y(t)\ge y_{\e}(t)$ for every $t\ge t_0$.
It is easy to check that
$$
y_{\e}(t) = y_0(t_0) \cosh c(t-t_0) + \frac{y_0'(t_0)-\e}{c} \sinh c(t-t_0) \,,
$$
for every $t,\e \in \RR$.
Hence
$$
y(t) \ge y_0(t_0) \cosh c(t-t_0) + \frac{y_0'(t_0)-\e}{c} \sinh c(t-t_0) \,,
$$
for every $t\ge t_0$ and $\e \in (0,y_0'(t_0))$.
If $\e \to 0$, we obtain
$$
y(t) \ge y_0(t_0) \cosh c(t-t_0) + \frac{y_0'(t_0)}{c} \sinh c(t-t_0) = y_0(t) \,,
$$
for every $t\ge t_0$.
This finishes the proof of the lemma.
\end{proof}

Lemma \ref{l:sturm} has the following direct consequence.

\begin{corolario}
\label{c:sturm}
Let us consider a positive constant $c$ and a function
$y$ satisfying
$y''(t)\ge c^2 y(t) >0$
and $y'(t_0) > 0$.
Then
$$
y(t) \ge y(t_0) \cosh c(t-t_0) \,,
$$
for every $t\ge t_0$.
\end{corolario}

\begin{proof}
Let us consider the function $y_0$ with
$$
\left\{ \begin{array}{ll}
y_0''(t)= c^2 y_0(t) \,,
& \;  \text{ if } \, t\ge t_0 \,,
\\
y_0(t_0) = y(t_0) \,,
& \;
\\
y_0'(t_0) = y'(t_0) \,.
& \;
\end{array}
\right.
$$
The first inequality in the following expression is obtained
by applying Lemma \ref{l:sturm} and
the first equality by solving the above differential equation:
$$
\aligned
y(t) & \ge y_0(t) = y_0(t_0) \cosh c(t-t_0) + \frac{y_0'(t_0)}{c} \sinh c(t-t_0)
\\
& \ge y_0(t_0) \cosh c(t-t_0)
= y(t_0) \cosh c(t-t_0) \,,
\endaligned
$$
for every $t\ge t_0$.
\end{proof}

The following result assures that if $K\le -c^2 < 0$,
there always exists a closed geodesic in every free homotopy class,
except for punctures,
in which is impossible to have one.

\begin{teorema}
\label{t:geodconknegativa}
Let us consider a Riemannian surface $S$,
which can be either simple bordered or without border.
Besides, $S$ must be complete and
with curvature $K\le -c^2 < 0$.
Fix an homotopically nontrivial closed curve $\a$ in $S$.
Then there exists a minimizing closed geodesic $\g \in [\a]$
if and only if $L([\a]) > 0$.
\end{teorema}

\begin{obs}
The conclusion of this Theorem
does not hold if we replace the hypothesis $K\le -c^2 < 0$ by the
weaker one $K < 0$, as shows the revolution surface of the graph of
$f(x)=1+e^x$ around the horizontal axis (with the standard metric
induced by the Euclidean metric in $\RR^3$).
\end{obs}

\begin{proof}
We deal first with nonbordered surfaces $S$.

If $L([\a]) = 0$,
it is clear that there does not exist a closed geodesic $\g \in [\a]$
with $L(\g) = L([\a]) = 0$, since $\a$
is an homotopically nontrivial closed curve in $S$.

Let us assume now that $L([\a]) > 0$.

Assume first that $S$ is not doubly connected.
Seeking for a contradiction, suppose that there not exist
a closed geodesic $\g \in [\a]$.
Then, by Proposition \ref{p:alternativa},
the curve $\a$ bounds a generalized puncture $E$ in $S$.
Since the curvature satisfies $K\le -c^2 < 0$,
this end $E$ is a Riemannian collared end,
by Theorem \ref{t:ends1}.

For each $r_0$ we define $g_{r_0}$ as the closed curve $\{r=r_0\}$.
It is easy to check that
$l(r):= L(g_r)= \int_0^{2\pi} G(r,\th) \, d\th$ satisfies
$l''(r) \ge c^2 l(r)$:
$$
\frac{\p^2 G}{\p r^2}\, (r,\th)
= - K(r,\th) G(r,\th)
\ge c^2 G(r,\th) > 0
$$
implies
$$
l''(r)
= \int_0^{2\pi} \frac{\p^2 G}{\p r^2}\, (r,\th) \, d\th
\ge \int_0^{2\pi} c^2 G(r,\th) \, d\th
= c^2 l(r) > 0 \,.
$$
Since $L([\a]) > 0$, there exist positive constants $c_0, r_1$
with $l(r) \ge c_0$ for every $r \ge r_1$. Hence, for every $r \ge r_1$,
$$
l'(r)
= l'(r_1) + \int_{r_1}^r  l''(t) \, dt
\ge l'(r_1) + \int_{r_1}^r  c^2 c_0 \, dt
= l'(r_1) + c^2 c_0 (r-r_1) \,,
$$
and consequently
$\lim_{r\to\infty} l'(r) = \infty$.
Since
$$
\lim_{r\to\infty} \int_0^{2\pi} \frac{\p G}{\p r}\, (r,\th) \, d\th
= \lim_{r\to\infty} l'(r) = \infty \,,
$$
there exist $r_2 \ge r_1$ and a set $A \subset [0,2\pi]$
with positive Lebesgue measure such that
$\p G/\p r (r_2,\th) > 0$ for every $\th \in A$.
Since $\p^2 G/\p r^2 (r,\th) > 0$,
the function $\p G/\p r (r,\th)$ increases in $r \ge r_2$
for each fixed $\th \in A$,
and consequently
$\p G/\p r (r,\th) \ge \p G/\p r (r_2,\th)> 0$ for every $\th \in A$ and $r \ge r_2$.
Hence, $G(r,\th)$ increases in $r \ge r_2$
for each fixed $\th \in A$.
By Corollary \ref{c:sturm},
$G(r,\th) \ge G(r_2,\th) \cosh c(r-r_2)$
for every $\th \in A$ and $r \ge r_2$.

Let us consider a curve $\s$ parametrized in the Riemannian collared end as
$\s(\th)=(r(\th),\th)$, with $\th \in [0,2\pi]$ and $r(\th) \ge R \ge r_2$.
Then
$$
\aligned
L(\s)
& = \int_0^{2\pi} \sqrt{r'(\th)^2 + G (r(\th),\th)^2} \; d\th
\ge \int_0^{2\pi}  G (r(\th),\th) \, d\th
\ge \int_A  G (r(\th),\th) \, d\th
\\
& \ge \int_A  G (R,\th) \, d\th
\ge \cosh c(R-r_2) \int_A  G (r_2,\th) \, d\th \,.
\endaligned
$$
Since $\int_A  G (r_2,\th) \, d\th$ is a positive constant independent of $R$,
there exists $r_3 > r_2$ with
$$
\cosh c(r_3-r_2) \int_A  G (r_2,\th) \, d\th > L(\a) \,.
$$
Hence, given any curve $\s$ parametrized in the Riemannian collared end as
$\s(\th)=(r(\th),\th)$, with $\th \in [0,2\pi]$ and $r(\th) \ge r_3$,
we have $L(\s) > L(\a)$.
Consequently, any curve $\s \in [\a]$ contained
in the region $\{ r \ge r_3 \}$ verifies $L(\s) > L(\a)$.
Then a minimizing sequence for $\a$ can not converge to $E$.
This fact contradicts Proposition \ref{p:alternativa}.

If $S$ is doubly connected,
the argument is similar except for the fact that
$\a$ bounds two collared ends.
Therefore, a minimizing sequence might not converge to an end in $S$;
but we can always extract a subsequence converging to some end in $S$.

Assume now that $S$ is simple bordered.
By Lemma \ref{l:complecion0},
$S$ is a subset of a complete Riemannian surface $R$.
The previous argument gives the desired result in $R$.
Then, Lemma \ref{l:complecion0} implies the result in $S$.
\end{proof}

The following lemma is a well known result,
but we include a direct proof by the sake of completeness.

\begin{lema}
\label{l:unicidad}
Let us consider a Riemannian surface $S$,
which can be either simple bordered or without border.
Besides, $S$ must be complete and
with curvature $K < 0$.
Then in each free homotopy class there exists at most
a closed geodesic, and if there exists, then it is minimizing.
Consequently, every generalized funnel is a funnel.
\end{lema}

\begin{proof}
By Lemma \ref{l:complecion0},
without loss of generality we can assume that $S$ is nonbordered.

Seeking for a contradiction, suppose that
there exist two freely homotopic closed geodesics $\g_1,\g_2$.

If $\g_1$ and $\g_2$ intersect at some point, they can not be
tangent at this point, since they are geodesics.
Then, $\g_1$ and $\g_2$ intersect at least at another point, since they are
freely homotopic.
Therefore, some segment of $\g_1$
and some segment of $\g_2$ determine a geodesic ``bigon" $B$
(a polygon with two sides)
with interior angles $\a,\b >0$.
Gauss-Bonnet Formula gives
$$
\iint_B K \, dA = \a + \b > 0 \,,
$$
which is a contradiction with $K < 0$.

Then $\g_1$ and $\g_2$ do not intersect. We consider the geodesic
segment $\s$ joining $x_1 \in \g_1$ with $x_2 \in \g_2$, which gives
the minimum distance between $\g_1$ and $\g_2$; then $\s$ meets
orthogonally to $\g_1$ and to $\g_2$. Let us consider a universal
covering map $\pi : \tilde S \longrightarrow S$. Fix a lift $\tilde
\g_1$ of $\g_1$ starting in $\tilde x_1$, a lift $\tilde \s$ of $\s$
starting in $\tilde x_1$, and finishing in $\tilde x_2$, and a lift
$\tilde \g_2$ of $\g_2$ starting in $\tilde x_2$. Then $\tilde \s$
meets orthogonally to $\tilde \g_1$ and to $\tilde \g_2$. If we
denote by $y_1$ and $y_2$, respectively, the endpoints of $\tilde
\g_1$ and $\tilde \g_2$, there exists a covering isometry $T :
\tilde S \longrightarrow \tilde S$ with $T(\tilde x_1)=y_1$ and
$T(\tilde x_2)=y_2$. We also have that with $T(\tilde \s)$ joins
$y_1$ and $y_2$, and meets orthogonally to $\tilde \g_1$ and to
$\tilde \g_2$. Consequently, $\tilde \g_1$, $\tilde \g_2$, $\tilde
\s$ and $T(\tilde \s)$ bound a geodesical quadrilateral $Q$ in
$\tilde S$ with four right angles. Gauss-Bonnet Formula gives
$$
- \iint_Q K \, dA = 2\pi - \frac{\pi}2 - \frac{\pi}2 - \frac{\pi}2 - \frac{\pi}2 = 0 \,,
$$
which is a contradiction with $K < 0$.

Then in each free homotopy class there exists at most
a closed geodesic.

Consider now a simple closed geodesic $\g$.
We need to prove that $L(\g)=L([\g])$.
Let us define $l:=L(\g)$.

The Riemannian metric can be expressed in
Fermi coordinates based on $\g$ as $ds^2= dr^2 + G(r,\th)^2 \, d\th^2$, where
$G(r,\th)$ satisfies
$$
\frac{\p^2 G}{\p r^2}\, (r,\th)
+ K(r,\th) G(r,\th)=0 \,,
\qquad G(0,\th)=1 \,,
\qquad \frac{\p G}{\p r}\, (0,\th)=0 \,.
$$
Since $\p^2 G / \p r^2 (r,\th) = - K(r,\th) G(r,\th) > 0$,
it follows that $G(r,\th)$ is a convex function on $r$ for each fixed $\th \in [0,l]$;
since $\p G / \p r (0,\th) = 0$, we deduce that for each fixed $\th \in [0,l]$,
$G(r,\th)$ attains its minimum value $1$ at $r=0$.

We prove now that any curve $\s \in [\g]$ verifies $L(\s) \ge L(\g)$.
Let us consider a fixed curve $\s \in [\g]$.
Without loss of generality we can assume that
$\s$ can be parametrized in Fermi coordinates as
$\s(\th)=(r(\th),\th)$, with $\th \in [0,l]$.
Then,
$$
L(\s) = \int_0^l \sqrt{r'(\th)^2 + G (r(\th),\th)^2} \; d\th
\ge \int_0^l  G (r(\th),\th) \, d\th
\ge \int_0^l  G (0,\th) \, d\th
= \int_0^l d\th = l = L(\g) \,.
$$
This shows that $L(\g)=L([\g])$ and hence $\g$ is minimizing.

In order to prove the last part of the lemma,
consider now a generalized funnel $F$ in $S$.
We have proved that there does not exist another closed geodesic
freely homotopic to the boundary of the generalized funnel.
Hence, the generalized funnel is a funnel.
\end{proof}

\begin{lema}
\label{l:punturasconknegativa}
Let us consider a Riemannian surface $S$,
which can be either simple bordered or without border.
Besides, $S$ must be complete and
with curvature $K\le -c^2 < 0$.
Then, every generalized puncture is a puncture.
\end{lema}

\begin{proof}
By Lemma \ref{l:complecion0},
without loss of generality we can assume that $S$ is nonbordered.

Seeking for a contradiction,
consider a generalized puncture which is not a puncture.
Then, its fundamental group is generated by a simple closed curve $\s$,
there is no minimizing simple closed geodesic $\g \in [\s]$,
and we have either:

$(i)$
$L([\s]) > 0$; then by Theorem \ref{t:geodconknegativa}
there exists a minimizing simple closed geodesic $\g \in [\s]$,
which is a contradiction,

\noindent or

$(ii)$
$L([\s]) = 0$ and
there exists a simple closed geodesic $\g \in [\s]$;
then by Lemma \ref{l:unicidad}
$\g$ is minimizing, which is a contradiction.
\end{proof}

\begin{teorema}
\label{t:disjuntas} Let us consider a complete Riemannian surface
$S$ and two disjoint nontrivial piecewise smooth simple closed
curves $a$ and $b$ in $S$, which are not freely homotopic. If
$\a\in [a]$ and $\b\in [b]$ are any choice of minimizing closed
geodesics, then they are disjoint as well.

Furthermore, if $\a\neq \b$ are freely homotopic minimizing simple
closed geodesics in $S$, then they are disjoint.
\end{teorema}

\begin{obs}
We have some examples that show that the conclusion of the previous
Theorem does not hold if either $\a$ or $\b$ are not minimizing
geodesics. See for example the figure below:
$$
\psfrag{a}{$\a$} %
\psfrag{b}{$\b$} %
\epsfysize=5cm\epsffile{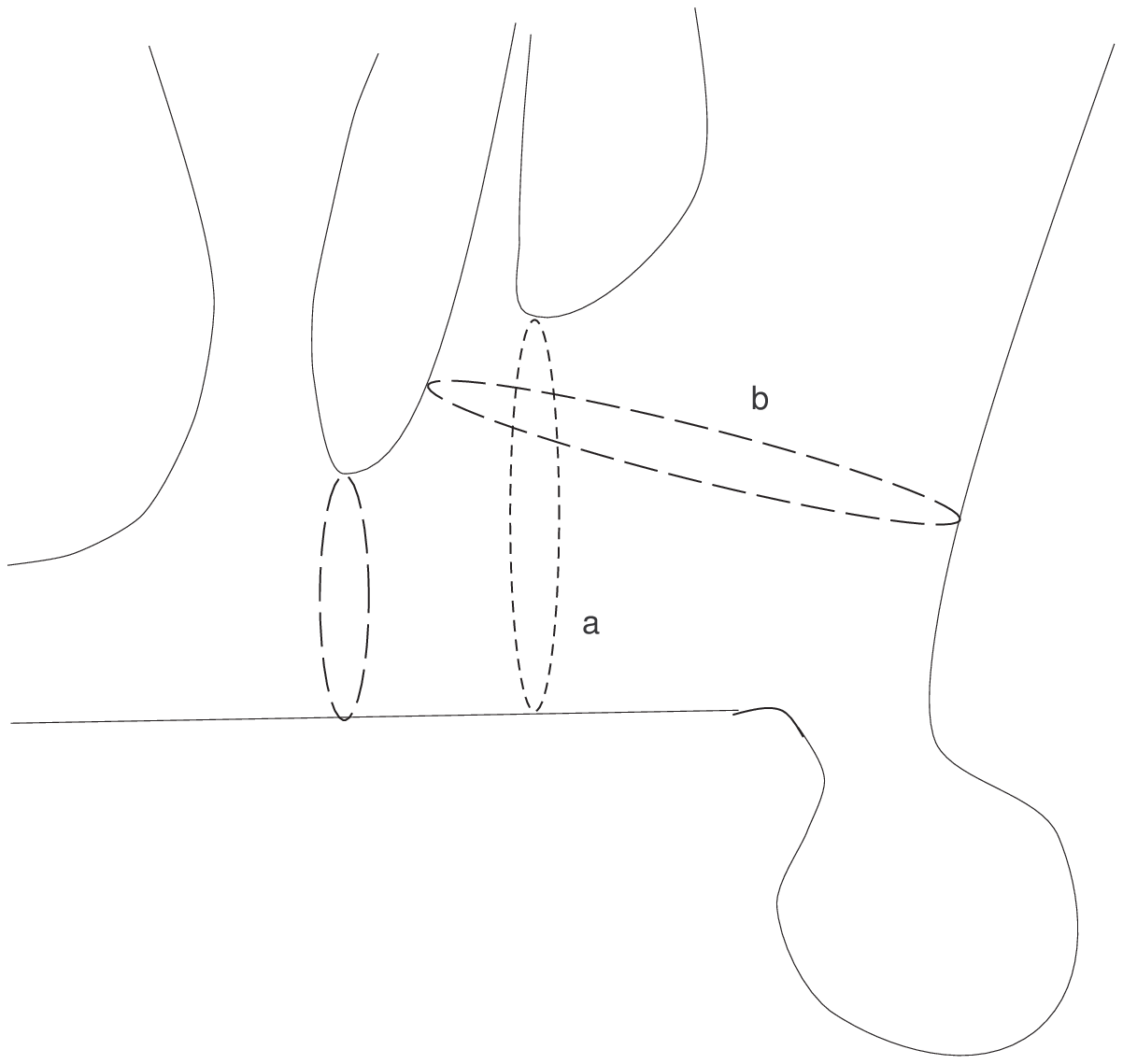}
$$
\end{obs}

\begin{proof}
First, let us assume that the curves $a$ and $b$ are not freely homotopic.
Without loss of generality we can assume that
$a$ and $b$ are disjoint nontrivial smooth simple closed curves in $S$,
since in other case we can modify them slightly in order to obtain all these facts.
Since $\a\in [a]$ and $S$ is a surface, by Baer's
Theorem (see e.g. \cite{St} or \cite{PS}) there exists an isotopy, that is to say, a
continuous family of diffeomorphisms $f_t:S\longrightarrow S$, such that
$f_0$ is the identity, and $f_1(a)=\a$. Let us define
$b_1:=f_1(b)$, which is a simple closed curve freely homotopic to
$b$. As $f_1$ is bijective and $a$ and $b$ are disjoint curves,
then $b_1$ and $\a$ are disjoint too.

Seeking for a contradiction, let us assume that $\a\cap\b \neq
\varnothing$. If they do intersect each other tangentially then they
should coincide, and this is not possible since they are not freely
homotopic. Therefore, they must intersect transversally. Since
$b_1\in[\b]$, there exists an smooth homotopy $F:A\longrightarrow
S$, where $A$ is the annulus $\{z \in \CC: \, 1\le |z| \le 2\}$,
such that $F(e^{i\theta})= b_1(\theta)$ and $F(2e^{i\theta})=
\b(\theta), \text{ with } \theta \in [0,2\pi]$.

No connected component $\g$ of $F^{-1}(\a)$ can be a closed curve in
$A$, since it should be either trivial or freely homotopic to
$F^{-1}(\b)$. In this case, $F(\g)=\a$ would also be either trivial
or homotopic to $F(F^{-1}(\b))=\b$, and this is contradiction with
our hypothesis. Therefore, $F^{-1}(\a)$ must contain an arc $\sigma$
joining two points $z_1$ and $z_2$ in $\{z\in \CC: \,
|z|=2\}=F^{-1}(\b)$. As $A$ is a planar domain, one of the two arcs
joining $z_1$ and $z_2$ in $F^{-1}(\b)$ (let us denote such arc by
$\eta$), is homotopic to $\sigma$. This fact implies that the
geodesics $\a$ and $\b$ intersect in $F(z_1)$ and $F(z_2)$ and that
$F(\sigma)$ and $F(\eta)$ are homotopic. Since $\a$ and $\b$ are
minimizing, $L(F(\sigma))=L(F(\eta))$ and so, from $\a$
we can construct a new curve $\tilde{\a}\in[\a]$ with the same
length by replacing the arc $F(\sigma)$ by the arc $F(\eta)$. By
means of a smooth modification in small neighborhoods of $F(z_1)$
and $F(z_2)$ we can obtain a shorter curve freely homotopic to $\a$,
which is contradiction with the fact that $\a$ is minimizing.

Now, we will deal with the second part of the Theorem. Once
again we are going to seek for a contradiction: let us assume that
$\a\cap\b \neq \varnothing$. If they do intersect in a single
point, as they are in the same homotopy class, they must
intersect each other tangentially and therefore $\a=\b$. If
they do intersect in several points, then they must intersect
transversally. The argument in the previous
case allows to obtain a shorter curve freely homotopic to
$\a$, which is contradiction with the fact that $\a$ is
minimizing.
\end{proof}

\section{The main results.}

\begin{prop}
\label{p:domgeod}
Every geodesic domain in any complete orientable Riemannian surface
is a finite union (with pairwise disjoint interiors) of generalized Y-pieces.
\end{prop}

\begin{obs}
\label{o:domgeod}
The argument in the proof of Proposition \ref{p:domgeod}
also proves the following:

Every complete orientable Riemannian surface without generalized
funnels and with finitely generated fundamental group, which is
neither simply nor doubly connected nor homeomorphic to a torus, is
a finite union (with pairwise disjoint interiors) of generalized
Y-pieces.
\end{obs}

\begin{proof}
Let us fix a geodesic domain $G$
in a complete orientable Riemannian surface $S$.
We denote by $\g_1, \g_2, \dots, \g_k$
the minimizing simple closed geodesics in $\partial G$.
Since $G$ is a simple complete bordered Riemannian surface,
by Lemma \ref{l:complecion0}
it is a subset of a complete Riemannian surface $R$.

In particular, $R$ is a complete orientable topological surface.
Since $R$ contains a geodesic domain, $R$ is
neither simply nor doubly connected nor homeomorphic to a torus.
Then, by Theorem \ref{t:descomposicion02},
$R$ is the union of topological Y-pieces $\{Y_n\}$ and cylinders $\{C_n\}$.
The fundamental group of $R$ is isomorphic to the
fundamental group of $G$, and therefore it is finitely generated;
then there are only a finite number of
topological Y-pieces and cylinders.
We denote by $\{\eta_m\} \subset R$ the set of
pairwise disjoint simple closed curves in $\cup_{n} \p Y_n$.
Without loss of generality we can assume that the curves are numbered
such that $\eta_j \in [\g_j]$ for each $1 \le j \le k$.

We want to change the curves
$\eta_j$ by minimizing simple closed geodesics.
For each $1 \le m \le k$, we replace $\eta_m$ by $\g_m$
(Lemma \ref{l:complecion0} gives that $\g_m$ are also
minimizing simple closed geodesic in $R$).
For each $m > k$,
let us choose a minimizing simple closed geodesic $\g_m \in [\eta_m]$, if it exists.
In other case, by Proposition \ref{p:alternativa},
the curve $\eta_m$ bounds a generalized puncture in $S$
and we define $\g_m := \varnothing$.
By Theorem \ref{t:disjuntas}, the minimizing simple closed geodesics
$\{\g_m\} \subset G$ are pairwise disjoint;
then, they split $G$ in the required
finite union of generalized Y-pieces
(if for some $m$ we have $\g_m = \varnothing$, the corresponding $Y$-piece in $R$
is a generalized $Y$-piece in $G$).
\end{proof}

The following theorem is the main result of this paper.
It generalizes an already known result for constant negatively curved surfaces
to arbitrary surfaces with no restricition of curvature at all.

\begin{teorema}
\label{t:descomposicion}%
Every complete orientable Riemannian surface which is neither
simply nor doubly connected nor homeomorphic to a torus is the
union (with pairwise disjoint interiors) of generalized Y-pieces,
generalized funnels and halfplanes.
\end{teorema}

\begin{obs}
\label{r:descomposicion}%
If there are several freely homotopic minimizing simple closed
geodesics which bound a generalized funnel,
we will see in the proof that any of them can be chosen
as border of this generalized funnel.
\end{obs}

\begin{proof}
We assume first that the fundamental group of $S$
is finitely generated.
If $S$ has not generalized funnels, then Remark \ref{o:domgeod} gives the result.
If $S$ has generalized funnels $\{F_j\},$
then the closure of $S\setminus \cup_j F_j$ is a geodesic domain;
Proposition \ref{p:domgeod} gives the result in this case.

Let us consider a surface $S$
with infinitely generated fundamental group,
and fix a point $p\in S$. Next,
we will take an increasing sequence of positive numbers
$\{r_n\}$ so that $\lim_{n\to\infty}r_n=\infty$. For each $r_n$
we intend to associate a geodesic domain $G_n$ to the ball
$B(r_n)$ centered in $p$ with radius $r_n$.

The boundary of $B(r)$ is a finite union of pairwise disjoint
simple closed curves except for $r\in A$ with $A$ a countable set.
Since $S$ is of infinite type, we can always find a positive
number $r_1\notin A$ such that the fundamental group of the ball
$B(r_1)$ has, at least, two generators. We choose $r_n \notin A$
with $r_n>\max\{r_{n-1}, n\}$. As $r_n>r_1$, the fundamental group
of $B(r_n)$ has, at least, two generators as well. Since
$r_n\notin A$, the boundary of $B(r_n)$ is a finite union of
pairwise disjoint simple closed curves $\{\eta_i^n\}_{i\in I_n}$.
In order to construct its geodesic domain $G_n$, our goal is to
relate a minimizing geodesic $\g_i^n$ to each curve
$\eta_i^n\subseteq\partial B(r_n)$, and we do it inductively as it
follows. There are two possibilities:
\begin{enumerate}
\item There not exists any minimizing simple closed geodesic in $[\eta_i^n]$.
In this case $\g_i^n:=\varnothing$. %
\item There exists at least one minimizing simple closed geodesic
in $[\eta_i^n]$. If $n>1$ and there is $j\in I_{n-1}$ such that
$\g_j^{n-1}\in[\eta_i^n]$, then $\g_i^n:=\g_j^{n-1}$. Otherwise
(notice that this situation includes the case $n=1$), choose
$\g_i^n$ as any of the minimizing simple closed geodesics in
$[\eta_i^n]$.
\end{enumerate}

$G_n$ is the geodesic domain limited by all these geodesics
$\{\g_i^n\}_{i\in I_n}$. By construction, $G_n \subseteq G_{n+1}$.

\spb

Before going on with the proof, we need the following lemma:

\begin{lema}
\label{l:fonilgen} %
If there exists some positive number $N$ such that $\g$ is a
minimizing simple closed geodesic contained in $\partial G_n$ for
every $n>N$, then $\g$ is the border of a generalized funnel in $S$.
\end{lema}

\begin{proof}
For $n> N$, let us consider the simple closed curve $\eta_n
\subseteq \partial B(r_n)$ which is freely homotopic to $\g$.
Since $\lim_{n\to
\infty}\dist(p,\eta_n)=\lim_{n\to\infty}r_n=\infty$, and $\eta_n$
belongs to a single nontrivial freely homotopy class for every
$n>N$, Theorem \ref{t:ends2} gives that $\{\eta_n\}$ converges to
a collared end $F$; since its border $\g$ is a minimizing simple
closed geodesic, $F$ must be a generalized funnel.
\end{proof}

Now, let us continue with the proof of Theorem
\ref{t:descomposicion}. We can take a subsequence of radii (by
simplicity of the notation we will denote this subsequence just
like the whole sequence) such that $G_n \subset G_{n+1}$ and
besides, $\partial G_n \cap \partial G_{n+1}$ is either the empty
set or a union of minimizing simple closed geodesics, each of them
is the border of a generalized funnel.

Let us define $H_n$ as the closed set obtained as the union of
$G_n$ and the generalized funnels whose borders belong to
$\partial G_n$, and $H:=\cup_n H_n$. Notice that due to the
properties of $G_n$, each $H_n$ is contained in the interior of
$H_{n+1}$. If $S=H$, then there is nothing else to prove.
Otherwise, $S\setminus H$ is a closed non-empty set.

By Proposition \ref{p:domgeod}, each connected component of
the closure of $G_{n+1}\setminus G_n$ is a finite
union (with pairwise disjoint interiors) of generalized Y-pieces.
Therefore,
in order to finish the proof, we just have to see that every connected
component $J$ of $S\setminus H$ is a halfplane, that is to say: $J$
is a simply connected set and $\partial J \subseteq \partial H$
is a unique nonclosed simple geodesic. From now on, by
simplicity in the notation and as there is no possible confusion, we
will denote $\g^n_i\subset
\partial H_n$ by $\g_n$.

Next, we state a lemma that we will need along the proof:

\begin{lema}
\label{l:gsc} %
Let $\sigma$ be a nontrivial simple closed curve in $S$. If there
is a minimizing simple closed geodesic $\g\in[\sigma]$ contained
in $B(r_n)$, then either $\g$ is contained in $G_n$ or it is
freely homotopic to some geodesic in $\partial G_n$.
\end{lema}

\begin{proof}
Let us assume that $\g$ is not contained in $G_n$. Then, there are
two possibilities: either $\g \cap G_n = \varnothing$ or $\g \cap
\partial G_n \neq \varnothing$. In the first case, $\g$ is contained
in a doubly connected set whose borders are $\g_n\subset
\partial G_n$ and $\eta_n\subset
\partial B(r_n)$ and therefore $\g$ is freely homotopic to both of them.

We finish the proof by showing that $\g$ cannot intersect
$\partial G_n$. Seeking for a contradiction, assume that $\g\cap
\partial G_n \neq \varnothing$; then, $\g\cap\g_n\neq\varnothing$
for some minimizing closed geodesic $\g_n\subset\partial G_n$.
Since $\g$ and $\g_n$ are minimizing simple closed geodesics,
Theorem \ref{t:disjuntas} implies that $\g\cap\eta_n\neq\varnothing$,
where $\eta_n$ is the closed curve in $\partial B(r_n)$ with
$\g_n\in[\eta_n]$. This is the required contradiction, since
$\g\subset B(r_n)$.
\end{proof}

Going on with proof of Theorem \ref{t:descomposicion} we will
see that $J$ is simply connected.
In order to do so, let us prove that
every simple closed curve contained in $J$ must be trivial, since
every topological obstacle must be in $H$: Let us consider a
nontrivial simple closed curve $\sigma$ in $J$. By Proposition
\ref{p:alternativa}, there are two possibilities: %
\begin{enumerate}
\item There exists a minimizing simple closed geodesic
$\g\in[\sigma]$. Consider $r_n$ with $\g\subset B(r_n)$; we prove
now that $\g\subset H_{n+1}$. If $\g$ is not contained in $G_n$,
then by Lemma \ref{l:gsc} it is freely homotopic to some geodesic
in $\partial G_n$. If $\g$ bounds a generalized funnel, then
$\g\subset H_n$. If a curve in $\partial G_n$ is freely homotopic
to some curve in $\partial G_{n+1}$ then it must bound a
generalized funnel. Since $\g$ does not bound a generalized
funnel, it is not freely homotopic to any geodesic in $\partial
G_{n+1}$; hence, Lemma \ref{l:gsc} gives $\g\subset G_{n+1}$, since
$\g\subset B(r_n)\subset B(r_{n+1})$. Hence, in any case,
$\g\subset H_{n+1}$, and it is not freely homotopic to any curve
in $\partial H_{n+1}$. Therefore, $\sigma$ must intersect
$H_{n+1}$, which is a contradiction with $\sigma\subset
J$.%
\item The curve $\sigma$ bounds a generalized puncture. Then,
there exists some neighborhood of this collared end contained in
some $H_n$. Since $\sigma$ is not freely homotopic to any geodesic
in $\partial H_n$, $\sigma$ must intersect $H_n$, which is again a
contradiction with $\sigma\subset J$.
\end{enumerate}

Next, we will prove that $\partial J$ is a geodesic. Let us fix a
point $q \in \partial J$; we want to prove that $q$ belongs to a
geodesic arc $\g$ contained in $\partial J$: There exist points
$q_n \in \g_n\subseteq\partial H_n$ converging to $q$. Let us
consider now the sequence $\{v_n\}$ of tangent vectors to $\g_n$
in $q_n$. Notice that this latest sequence must converge to a
certain vector $v$, since otherwise the geodesics $\g_n$ would
intersect. As geodesics are solutions of a system of ordinary
differential equations and the initial data $\{q_n,v_n\}$ converge
to $\{q,v\}$, then the geodesics $\g_n$ converge uniformly in some
neighborhood $U$ of $q$ to the geodesic $\g$ whose tangent vector
in $q$ is $v$.

Now, let us prove that $\g\cap U$ is contained in $\partial J$.
Choose a point $q' \in \g\cap U$. On the one hand, $q' \notin
\text{ext}H$, since there exists a sequence of points in
$\g_n\subset H_n$ converging to $q'$. On the other hand, $q'$ does
not belong to $H$ either, since if it did, there would exist some
$H_m$ containing $q'$ for some positive integer $m$, and this is a
contradiction with $H_n$ contained in the interior of $H_{n+1}$
for every $n$.

From the previous argument we also deduce that this geodesic can
be prolonged to infinity at both sides: if it had an endpoint $p$,
it is obvious that $p\in\partial J$, but as we have just seen, for
every point in the boundary of $J$ there exists a neighborhood $U$
such that $U\cap\partial J$ is a geodesic arc containing $p$.

In order to see that every connected component $\g \subseteq
\partial J$ is simple, we will prove that it is not closed and it does not intersect
itself transversally: If $\g$ were a simple closed geodesic, it
would be compact and as $\g_n$ locally uniformly converge to $\g$
then there would exist a positive integer $N$ and a collar $C$ for
$\g$ such that $\g_n \subset C$ for every $n\ge N$, and therefore
$\g_n\in [\g]$, which is a contradiction. If $\g$ intersected
itself transversally, there would exist some positive integer $N$
such that each $\g_n$ would intersect itself as well for every
$n\ge N$, and this is not possible since they are all simple. This
same argument also proves that two different geodesics contained
in $\partial H$ must be simple and pairwise disjoint.

To finish, there is only one fact to prove: $\partial J$ consists
of just one geodesic. Let us assume that there exist two simple
geodesics $\sigma_1, \sigma_2 \subset\partial J$. Let us consider
two points $q_1\in\sigma_1, q_2\in\sigma_2$, two simple connected
neighborhoods $V_1,V_2$ of $q_1$ and $q_2$ respectively, two
simple closed geodesics $\g_{n_1}\subset\partial H_{n_1},
\g_{n_2}\subset\partial H_{n_2}$, with $\g_{n_1}\cap
V_1\neq\varnothing, \g_{n_2}\cap V_2\neq\varnothing$ and $n_1 \neq
n_2$, and curves $\eta_1 \subset V_1, \eta_2\subset V_2$ joining,
respectively, $\g_{n_1}$ with $q_1$ and $q_2$ with $\g_{n_2}$. As
$J$ is path-connected, it is possible to construct the three
following curves: $\eta_3\subset J$ joining $q_1$ and $q_2$, $\eta
:= \eta_1 + \eta_3 + \eta_2$ and the closed curve $\beta := \eta +
\g_{n_2} - \eta + \g_{n_1}$.
Since $\b$ cannot bound a generalized
puncture, every closed
geodesic $\g \in [\beta]$ verifies $\g \cap \sigma_1 \neq
\varnothing$ and $\g \cap \sigma_2 \neq \varnothing$; in
particular, this means that every minimizing simple closed
geodesic do intersect $\partial J$. But, by Lemma \ref{l:gsc},
there must exist a minimizing closed geodesic in $[\b]$ entirely
contained in $H$, which is a contradiction.
\end{proof}

In fact, the proof of
Theorem \ref{t:descomposicion}
gives the following result.

\begin{teorema}
\label{t:descomposicionyexaucion}
Every complete orientable Riemannian surface
which is neither simply nor doubly connected nor homeomorphic to a torus
is the union (with pairwise disjoint interiors)
of generalized funnels, halfplanes
and a set $G$ which can be exhausted by geodesic domains.
\end{teorema}

The curvature of a Riemannian surface homeomorphic to a
torus can not verify $K<0$;
then, Theorem \ref{t:descomposicion},
Lemma \ref{l:unicidad} and
Lemma \ref{l:punturasconknegativa}
give directly the following result.

\begin{teorema}
\label{t:descomposicion2}
Every complete orientable Riemannian surface
with curvature $K\le -c^2 < 0$,
which is neither simply nor doubly connected
is the union (with pairwise disjoint interiors)
of generalized Y-pieces, funnels and halfplanes.
Furthermore, every generalized puncture is a puncture.
\end{teorema}

In order to deal with bordered surfaces,
we need a last definition.

\begin{definicion}
A \emph{finite cylinder}
is a bordered Riemannian surface
which is homeomorphic to $\SS^1 \times [0,1]$,
whose border is the union of two simple closed geodesics,
and at least one of them is minimizing.
\end{definicion}

\begin{teorema}
\label{t:descomposicion3}
Every simple complete orientable bordered Riemannian surface
which is neither simply nor doubly connected
is the union (with pairwise disjoint interiors)
of generalized Y-pieces, finite cylinders, generalized funnels and halfplanes.

Furthermore, there is a bijection between finite cylinders in the decomposition
and nonminimizing simple closed geodesics in the border.
\end{teorema}

\begin{proof}
Let $S$ be a simple complete orientable bordered Riemannian surface
which is not simply nor doubly connected, whose border is the union
of simple closed geodesics $\{\g_i\}_{i\in I}$.
By applying
Lemma \ref{l:complecion0}
we can construct another complete Riemannian
surface $R$ by gluing a neighborhood $F_i$ of a collared end to each $\g_i$,
such that $R=\cup_{i\in I} F_i \cup S$.
It is obvious that $R$ is not
homeomorphic to a torus, since it is not compact.

By Theorem \ref{t:descomposicion} we know that $R$ is the union
(with pairwise disjoint interiors) of generalized Y-pieces,
generalized funnels and halfplanes.
Furthermore, by Remark
\ref{r:descomposicion}, if $\g_i$ is a minimizing simple closed geodesic in $S$,
for some $i \in I$, we can choose the decomposition in such a way that $\g_i$
belongs to the border of a generalized funnel.

If $\g_i$ is a nonminimizing simple closed geodesic in $S$ for some $i\in I$,
Lemma \ref{l:complecion0} ($(7)$ and $(2)$)
guarantees both that there exists
the minimizing simple closed geodesic
$\g_i^0 \in [\g_i]$ and that it is contained in $S$.
Then, the funnel $F_i$ in this decomposition intersects $S$
in a finite cylinder whose border is $\g_i^0\cup \g_i$.

Consequently, we obtain the desired decomposition in $S$
if we restrict to $S$ this decomposition in $R$.
\end{proof}

\begin{teorema}
\label{t:descomposicion4}
Every simple complete orientable bordered Riemannian surface
with curvature $K\le -c^2 < 0$,
is the union (with pairwise disjoint interiors)
of generalized Y-pieces, funnels and halfplanes.
Furthermore, every generalized puncture is a puncture.
\end{teorema}

\begin{proof}
The proof follows the argument of the proof of Theorem
\ref{t:descomposicion3}, using Theorem \ref{t:descomposicion2},
instead of Theorem \ref{t:descomposicion}.

Lemma \ref{l:unicidad} gives that
in each free homotopy class there exists at most
a closed geodesic.
Consequently, there are not finite cylinders in the decomposition.

We only need to study
the simple complete orientable bordered Riemannian surfaces
with curvature $K\le -c^2 < 0$
which are simply or doubly connected.

Let $S$ be such a surface.
$S$ can not be simply connected:
Seeking for a contradiction, suppose that
$S$ is simply connected; then $\p S$ can be considered a geodesic triangle
with three angles equal to $\pi$, and Gauss-Bonnet Formula gives
$$
- \iint_S K \, dA = \pi - \pi - \pi - \pi = - 2 \pi \,,
$$
which is a contradiction with $K < 0$.

If $S$ is doubly connected, then
Lemma \ref{l:unicidad} gives that
$\p S$ is a single simple closed geodesic
and besides it is minimizing.
Then, $S$ is a funnel.
\end{proof}




\

\

\begin{tabular}{cccc}
\small
\parbox{7cm}{
Ana Portilla, Jos\'e M. Rodr\'{\i}guez, Eva Tour\'{\i}s\\
Departamento de Matem\'aticas\\
Escuela Polit\'ecnica Superior\\
Universidad Carlos III de Madrid \\
Avenida de la Universidad, 30\\
28911 Legan\'es (Madrid)\\
SPAIN\\} & &\hspace{2cm} & \small
\parbox{6cm}{
\quad\\
 \\
 \\
 \\
 \\
 \\
 \\
 \\}
\end{tabular}


\begin{thebibliography}{99}

\bibitem{AR} Alvarez, V., Rodr{\'\i}guez, J. M.,
Structure theorems for Riemann and topological surfaces,
{\it J. London Math. Soc.} {\bf 69} (2004), 153--168.


\bibitem{B} Busemann, H.,
{\it The Geometry of Geodesics.}
Academic Press, New York, 1955.


\bibitem{C} Chavel, I.,
{\it Eigenvalues in Riemannian Geometry.}
Academic Press, New York, 1984.

\bibitem{Ch} Chazarain, J.,
{\it Spectre des op\'erateurs elliptiques et flots hamiltoniens.}
In: S\'eminaire Bourbaki 1974-75, Lecture Notes in Mat. {\bf 514},
Springer, Berlin-Heidelberg, 1976.


\bibitem{CV} Colin de Verdiere, Y.,
Quasi-modos sur les vari\'et\'es Riemanniennes,
{\it Inv. Math.} {\bf 43} (1977), 15--52.

\bibitem{E} Eberlein, P.,
{\it Surfaces of nonpositive curvature.} Memoirs of the American
Mathematical Society 218, American Mathematical Society, Providence,
1979.

\bibitem{FM} Fern\'andez, J. L., Meli\'an, M. V.,
Escaping geodesics of Riemannian surfaces,
{\it Acta Math.} {\bf 187} (2001), 213--236.

\bibitem{G} Guillemin, V.,
Lectures on spectral theory of elliptic operators,
{\it Duke Math. J.} {\bf 44} (1977), 485--517.


\bibitem{M} Massey, W. M.,
{\it  Algebraic Topology: An Introduction.}
Harcourt, Brace and World, Inc., New York, 1967.



\bibitem{PS} Prasolov, V. V., Sossinsky,  A. B.,
{\it Knots, links, braids and 3-manifolds.}
AMS Translations of Mathematical Monographs 154, 1997.


\bibitem{St} Stillwell, J., {\it Classical topology and combinatorial group theory
(Second Edition).}
Springer-Verlag, New York, 1993.



\end{thebibliography}
\end{document}